\documentclass[11pt,reqno]{amsart}

\usepackage{amsmath}
\usepackage{amsxtra}
\usepackage{amscd}
\usepackage{amsthm}
\usepackage{amsfonts}
\usepackage{amssymb}
\usepackage{eucal}
\usepackage[all]{xy}
\numberwithin{equation}{section}
\newtheorem{theorem}{Theorem}[section]

\newtheorem{lemma}[theorem]{Lemma}
\newtheorem{corollary}[theorem]{Corollary}
\theoremstyle{definition}

\newcommand{\nc}{\newcommand}
\nc{\C}{{\mathbb C}}

\let\al\alpha

\newcommand{\Z}{{\mathbb Z}}

\newcommand{\Ref}[1]{{$($\ref{#1}$)$}}
\newcommand{\bean}{\begin{eqnarray}}
\newcommand{\eean}{\end{eqnarray}}
\newcommand{\be}{$$}
\newcommand{\bea}{\begin{eqnarray*}}
\newcommand{\eea}{\end{eqnarray*}}
\newcommand{\ee}{$$}
\def\ov[#1,#2]{\overset{\scriptstyle #1}{#2}}
\nc{\verm}{M_{k,l}}
\def\wsl{\widehat{\mathfrak{sl}}_2}
\nc{\slt}{{\widehat{\mathfrak{sl}}_2}}
\nc{\inte}{{L_{k,l}}}
\nc{\pkln}{{\mathcal P}_{k,l}^{(N)}}
\nc{\ckln}{{c_{k,l}^{(N)}}}
\def\qbin[#1;#2]{{\left[\matrix{\displaystyle #1}\\{\displaystyle #2}\endmatrix\right]}}
\def\geq{\ge}
\def\leq{\le}

\nc{\Vk}{{\mathfrak V}_k}

\def\xx[#1,#2]{{#1}^{(#2)}}
\def\yy[#1,#2,#3]{{#1}^{(#2)}_{#3}}

\renewcommand{\deg}{{\rm deg}\,}

\def\v(#1;#2){\Bigl({#1\atop#2}\Bigr)}

\nc{\alb}{{\boldsymbol{\alpha}}}
\nc{\beb}{{\boldsymbol{\beta}}}

\def\hookdownarrow%
{\Big\downarrow\kern-.14cm\smash{\raise.4cm\hbox{$\scriptstyle\cap$}}}

\def \wsw {\widehat{\mathfrak{sl}}_2}

\nc{\heis}{\widetilde{\mathfrak H}}
\nc{\dualW}{{W_k^{(M,N)}[l_1,l_2,l_3]}^*}
\nc{\dualWa}{{W_k[l_1,l_2,l_3]}^*}
\nc{\card}{{\#}\,}
\nc{\bs}{\boldsymbol}
\nc{\mc}{\mathcal}
\nc{\on}{\operatorname}
\nc{\mf}{\mathfrak}
\nc{\ol}{\overline}

\begin{document}

\title[Bosonic formulas]{Bosonic formulas for $\wsw$ coinvariants}

\author{M. Jimbo, T. Miwa and E. Mukhin}
\address{MJ: Graduate School of Mathematical Sciences,
The University of Tokyo, Tokyo 153-8914, Japan}\email{jimbomic@ms.u-tokyo.ac.jp}
\address{TM: Division of Mathematics, Graduate School of Scince,
Kyoto University, Kyoto 606-8502 Japan}\email{tetsuji@kusm.kyoto-u.ac.jp}
\address{EM: Dept. of Mathematics, University of
California, Berkeley, CA 94720}\email{mukhin@math.berkeley.edu}

\begin{abstract}
We derive bosonic-type formulas for the characters of $\wsw$ coinvariants.
\end{abstract}

\maketitle

\section{Introduction}
Let us fix a natural number $k$ which we call level.
Let $L_{l}$ be the integrable $\wsw$ module of level $k$ and of
weight $l$, generated by a highest weight vector $v_{l}$
such that
\be
h_0v_{l}=lv_{l},\quad f_0^{l+1}v_{l}=e_1^{k-l+1}v_{l}=f_{i}v_{l}=e_{i+1}v_{l}=h_iv_{l}=0,
\qquad i\in \Z_{<0},
\ee
where $e_i,f_i,h_i$ are the standard generators of $\wsw$.
Then the space of coinvariants $L_{k,l}^e$ is the quotient space
\be
L_{l}^e=L_{l}/(e_iL_{l},\;\;i\in\Z_{\geq 0}).
\ee
This quotient is naturally doubly graded by actions of $h_0$
and the degree operator $d$.
In this paper we present bosonic-type formulas for the corresponding
character
\be
\chi_{l}^e(q,z)=\on{Tr}_{L_{l}^e}q^dz^{h_0}.
\ee
The spaces of coinvariants $L_{l}^e$ were studied in the series of papers
\cite{FKLMM1}, \cite{FKLMM2}, \cite{FKLMM3} (the space
$L_{l}^e$ is denoted there by
$L_{k,l}^{(0,\infty)}$), where several other spaces with the same character
$\chi_{l}$ are given. One such space, a space of rigged
configurations, produced a fermionic-type formula for
$\chi_{l}^e(q,z)$, see Theorems 3.5.2, 3.5.3
in \cite{FKLMM3}. Another description via
``combinatorial paths'' led to a difference equation for
$\chi_{l}^e(q,z)$. This difference equation was the main theme in
all proofs of the papers \cite{FKLMM1}, \cite{FKLMM2},
\cite{FKLMM3}.
It also plays the central role in this paper.
Let us recall the construction.

For $0\leq i\leq l\leq k$, denote by $\mc C_{l}[i]$ the set of pairs
$({\bf a};{\bf b})=(a_0,a_1\dots;b_0,b_1,\dots)$, where
$a_i,b_i$ are non-negative integers
such that only finitely many are different from zero, and
\bean\label{planes}
a_0=i,\quad a_r+b_{r+1}+a_{r+1}\leq k,\quad b_r+a_r+b_{r+1}\leq k,  \quad
\sum_{s=m}^n b_s \leq k+\sum_{s=m+1}^{n-2}a_s,
\eean
where  $r\geq 0$, $-1\leq m<n\leq \infty$. Here
we set $b_\infty=k$, $b_{-1}=l$, $b_0=a_{-1}=a_{\infty}=0$.
Then by Corollary 5.4.10  of \cite{FKLMM2},
\be
\chi_{l}(q,z)=\sum_{i=0}^l \chi_{i,l}(q,z,z^{-1}),\qquad
\chi_{i,l}(q,z_1,z_2):=\sum_{({\bf a},{\bf b})\in\mc C_{l}[i]}q^{\sum
j(a_j+b_j)}z_1^{\sum b_j}z_2^{\sum a_j+b_j}.
\ee

By Proposition 3.3.1 in \cite{FKLMM2},
we have a difference equation of the form
\be
\chi_{i,l}(q,z_1,z_2)=\sum_{i',l'}
M_{i,l}^{i',l'}(q,z_1,z_2)\chi_{i',l'}(q,z_1,qz_2), \qquad
\chi_{i,l}(q,z_1,0)=\delta_{i,0}\delta_{l,0},
\ee
where the matrix $M$ is given by \Ref{matrix} below.

The set $\mc C_{l}[i]$ can be thought of as the set of integer points
in a polytope in an infinite-dimensional space cut out by the
hyperplanes \Ref{planes}.
Then one may expect an existence of a bosonic-type formula
for $\chi_{i,l}$ in the spirit of \cite{KhP}, 
written as a sum over
vertices of the characters of the corresponding cones.
Such a bosonic formula was studied in a simpler situation 
in \cite{FL}. 
However, our polytope 
is rather complicated 
and a direct approach does not look promising.

Instead we use the difference equation
to reduce our infinite-dimensional problem
to a problem in two dimensions. 
The non-zero entries of the matrix $M$ are written as 
$M_{i,l}^{i',l'}=m_1^im_2^lm_3^{i'}m_4^{l'}$,
where the $m_j$ are monomials of the form 
\bean\label{mon form}
q^{\al_j}z_1^{\beta_j}z_2^{\gamma_j},\qquad\al_j,\beta_j,\gamma_j\in\Z.
\eean
Moreover, for a generic $(i,l)$, the set of
$(i',l')$ such that $M_{i,l}^{i',l'}\neq 0$ is a pentagon
(see Figure \ref{pentagon}).
Given a monomial $n_1,n_2$ of the form \Ref{mon form}, 
these properties allow us to rewrite 
$\sum_{i',l'}M_{i,l}^{i',l'} n_1^{i'}n_2^{l'}$
as a sum of 5 rational functions 
corresponding to the vertices of the pentagon.
We denote
these five terms by $A(n_1^{i}n_2^{l}),\cdots, E(n_1^{i}n_2^{l})$.
Then we can formally write
\begin{eqnarray}
\chi_{i,l}=\lim_{N\rightarrow \infty}
(A+B+C+D+E)^N (\delta_{i,0}\delta_{l,0}).
\label{formal}
\end{eqnarray}
Expanding the RHS we obtain an expression for the character as a  
sum over all possible non-commutative 
monomials in the operators $A,\cdots,E$. 

This na\i ve expansion suffers from difficulties for two reasons. 
First, the number of terms grows exponentially as $N\rightarrow\infty$. 
Second, there are terms which contain zero denominators. 
The aim of the present paper is to find a formula for $\chi_{i,l}$ 
as a sum over a subset of monomials, such that 
the number of terms grows polynomially and that all terms are well defined. 
Informally the existence of such a formula 
means that a huge cancellation takes place.  
The result is given in Theorem \ref{main}. 
In the preceding sections we explain the origin of this formula. 
For that purpose, an appropriate language is provided by what we call 
summation graph; see Figure \ref{add graph}.
Monomials in the operators 
$A,B,C,D,E$ are identified with directed paths in this graph, 
and the cancellation pattern can be visualized. 
Terms corresponding to the same path (of infinite length) 
can be summed up with the use of Jackson's 
${}_6\Phi_5$ formula (see \Ref{id} below). 
Then, all terms which cancel do so in pairs. 
Therefore, we need to perform 
no additional summation to observe the desired cancellation. 
We emphasize that the actual proof of Theorem \ref{main} 
is done by a direct computation and is logically 
independent of these considerations. 

Our final answer for $\chi_{i,l}$ given in Section \ref{expl section}
is a sum of 18 families of meromorphic functions.  
The functions in each family 
are parametrized by 3 non-negative integers.
Each of them has a factorized form 
\be
m_0m_1^im_2^lm_3^k\frac{\prod_r(1-f_r)}{\prod_s(1-g_s)},
\ee
where $m_j,f_r,g_s$ are monomials as in \Ref{mon form}.
Note that the structure of this formula does not depend on $i,l$
in contrast to the fermionic formulas.

We expect that one can write in a similar fashion bosonic
formulas for solutions of certain class of difference equations.
Another example of such equations and the corresponding formulas are
given in \cite{FJLMM}.

{\bf Acknowledgments.} The authors would like to thank B. Feigin for
many stimulating discussions.
This work is partially supported by
the Grant-in-Aid for Scientific Research (B)
no.12440039, Japan Society for the Promotion of Science.

\section{The difference equation}
Fix a natural number $k\in\Z_{>0}$. Define a
matrix
$M(q,z_1,z_2)=(M_{i,l}^{i',l'})$ where $0\leq i\leq l\leq k$, $0\leq
i'\leq l'\leq k$, by the formulas
\bean\label{matrix}
M_{i,l}^{i',l'}=\begin{cases}(qz_1z_2)^{l'-i'}z_2^i &
{\rm if}\;\; l-i\leq i'\leq l'\leq k-i;\cr
(qz_1z_2)^{l'-l+i}z_2^i & {\rm if}\;\; i'<l-i\leq l'\leq k-i;\cr
0&{\rm otherwise}. \end{cases}
\eean
It is a square matrix
of size $(k+1)(k+2)/2$.
We denote $M(q,z_1,0)$ simply by $M(0)$.
All entries of the matrix $M(0)$ are equal to $0$ or $1$.

{}For example, for $k=1$,
we have
\be
M=
\left(\begin{matrix}
1 & qz_1z_2 & 1\\
0 & 1 & 1\\
z_2 & 0 & 0
\end{matrix} \right), \qquad
M(0)=
\left(\begin{matrix}
1 & 0 & 1\\
0 & 1 & 1\\
0 & 0 & 0
\end{matrix} \right),
\ee
where the components are arranged in the order
$(i,l)=(0,0),(0,1),(1,1)$.
Consider the system of difference equations
\bean\label{diff eqn}
\chi(q,z_1,z_2)=M(q,z_1,z_2)\chi(q,z_1,qz_2)
\eean
for a vector $\chi(q,z_1,z_2)=(\chi_{i,l})$ where $0\leq i\leq
l\leq k$, and $\chi_{i,l}$ is a formal power series
in the variables $q,z_2$, whose coefficients are formal 
Laurent series in $z_1^{-1}$. 
We impose the following initial condition
\bean\label{init cond}
\chi_{i,l}(q,z_1,0)=\delta_{i,0}\delta_{l,0}.
\eean

\begin{lemma}\label{existence}
The system \Ref{diff eqn}, \Ref{init cond} has a unique solution
in $\C[[q,z_2]]((z_1^{-1}))^m$ where
$m=(k+1)(k+2)/2$.
\end{lemma}
\begin{proof}
Let us write 
$\chi(q,z_1,z_2)=\sum_{j=0}^\infty f_j(q,z_1)z_2^j$, 
where $f_j(q,z_1)\in \C[[q]]((z_1^{-1}))^m$. 
By \Ref{init cond}, $f_0(q,z_1)_{i,l}=(\delta_{i,0}\delta_{l,0})$.
Comparing the coefficients of $z^n$ in \Ref{diff eqn}, we obtain 
\be
f_n(q,z_1)=q^n M(0)f_n(q,z_1)+\dots,
\ee
where the dots denote terms which depend on $f_j(q,z_1)$ with $j<n$.
There exists the inverse matrix $(\on{Id}-q^nM(0))^{-1}$
whose coefficients are
power series in $q$ alone, 
and therefore the vector $f_n(q,z_1)$ is uniquely determined 
via $f_j(q,z_1)$ with $j<n$.
\end{proof}

{}For $N\in\Z_{\geq 0}$, define vector valued functions
$\chi^{(N)}(q,z_1,z_2)=(\chi^{(N)}_{i,l})$ ($0\leq i\leq l\leq k$)
recursively by
\bean\label{recursion}
\chi^{(N+1)}(q,z_1,z_2)&=&M(q,z_1,z_2)\chi^{(N)}(q,z_1,qz_2),\\
\chi^{(0)}_{i,l}(q,z_1,z_2)&=&\delta_{i,0}\delta_{l,0}.\notag
\eean
The function $\chi^{(N)}_{i,l}$ is a polynomial in 
$q,z_1,z_2$ with 
non-negative integer coefficients.

\begin{lemma}
As $N$ tends to infinity,
the limit of $\chi_{i,l}^{(N)}(q,z_1,z_2)$ exists 
in $\C[[q]][z_1,z_2]$ and is equal to $\chi_{i,l}(q,z_1,z_2)$.
\end{lemma}
\begin{proof}
Note that
\be
M(0)\chi^{(0)}(q,z_1,z_2)=\chi^{(0)}(q,z_1,z_2),
\ee
and that for $m\geq n$, $M(q,z_1,q^mz_2)\equiv M(0)\,\bmod\,q^n$. Therefore,
for $N\geq n$ we have,
\begin{alignat*}{2}
\chi^{(N)}(q,z_1,z_2)
&=
M(q,z_1,z_2)M(q,z_1,qz_2)\dots M(q,z_1,q^{N-1}z_2)\chi^{(0)}(q,z_1,q^Nz_2)
\\
&\equiv
M(q,z_1,z_2)M(q,z_1,qz_2)\dots
M(q,z_1,q^{n-1}z_2)M(0)^{N-n}\chi^{(0)}(q,z_1,z_2)
\\
&\equiv
M(q,z_1,z_2)M(q,z_1,qz_2)\dots
M(q,z_1,q^{n-1}z_2)\chi^{(0)}(q,z_1,z_2)
\\
&\phantom{M(q,z_1,z_2)M(q,z_1,qz_2)\dots
M(q,z_1,q^{n-1}z_2)\chi^{(0)}(q,z_1,z_2)}
\bmod q^n.
\end{alignat*}
Therefore
the limit of $\chi_{i,l}^{(N)}(q,z_1,z_2)$ exists as $N$ tends to infinity.
Then it is equal to $\chi_{i,l}(q,z_1,z_2)$ by the uniqueness 
part of Lemma \ref{existence}.
\end{proof}

\begin{corollary}\label{integercoeff}
The components of $\chi(q,z_1,z_2)$ are power series in $q$ with
coefficients in $\Z_{\geq 0}[z_1,z_2]$.
\end{corollary}

The main purpose of this paper is to derive a bosonic-type formula for
$\chi(q,z_1,z_2)$.

\section{Bosonic formula for a triangle and a rectangle}
\begin{figure}
\begin{picture}(30,120)(40,-10)
\put(0,0){\vector(1,0){150}}
\put(0,0){\vector(0,1){100}}
\put(0,0){\line(1,1){100}}

\put(0,43){\line(1,0){46}}
\put(0,46){\line(1,0){46}}
\put(0,86){\line(1,0){86}}

\put(46,43){\line(0,1){3}}
\put(86,43){\line(0,1){3}}

\thicklines

\put(46,46){\line(1,1){40}}
\put(46,46){\line(1,0){40}}
\put(86,46){\line(0,1){40}}

\put(46,43){\line(0,-1){43}}
\put(46,0){\line(1,0){40}}
\put(86,0){\line(0,1){43}}
\put(46,43){\line(1,0){40}}

\put(-10,95){$i'$}
\put(150,-10){$l'$}
\put(-10,0){$0$}

\put(-23,49){$l-i$}
\put(-25,83){$k-i$}
\put(-41,36){$l-i-1$}
\put(40,-10){$l-i$}
\put(80,-10){$k-i$}

\put(36,3){$A$}
\put(89,3){$C$}
\put(34,50){$B_1$}
\put(89,48){$D_1$}
\put(48,32){$B_2$}
\put(89,38){$D_2$}
\put(89,78){$E$}
\end{picture}
\caption{Region of summation in \Ref{diff eqn},
\Ref{recursion}}\label{pentagon}
\end{figure}
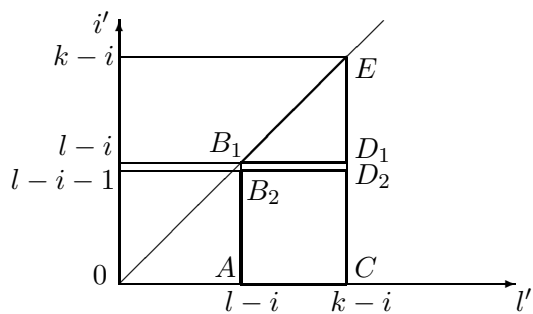

The sum in the difference equations \Ref{diff eqn}, \Ref{recursion}
is taken over the union of a triangle $B_1D_1E$ and a
rectangle $AB_2D_2C$ shown in Figure \ref{pentagon}.
We divide the pentagon into two regions
(the triangle $B_1D_1E$ and the rectangle $AB_2D_2C$)
since the formulas we apply will be different for these two.

The following two lemmas give bosonic type expressions for a sum over
integer points inside a triangle and a rectangle.
These lemmas are easily proved by a direct computation.
\begin{lemma}\label{tri id lemma}
Let $a,b$ be integers such that $a\leq b$. Then
\bean\label{tri id}
\sum_{a\leq n\leq m\leq b} x^my^n= \frac{(xy)^a}{(1-x)(1-xy)}+
\frac{x^by^a}{(1-x^{-1})(1-y)}+\frac{(xy)^b}{(1-y^{-1})(1-(xy)^{-1})}.
\eean
\end{lemma}

\begin{lemma}\label{rect id lemma}
Let $a,b,c,d$ be integers such that $a\leq b+1$ and $c\leq d+1$.
Then
\bean\label{rect id}
\sum_{{a\leq m\leq b}\atop{c\leq n\leq d}}
x^my^n&=&\frac{x^ay^c}{(1-x)(1-y)}+\frac{x^ay^d}{(1-x)(1-y^{-1})}
\notag\\
&+&\frac{x^by^c}{(1-x^{-1})(1-y)}+\frac{x^by^d}{(1-x^{-1})(1-y^{-1})}.
\eean
\end{lemma}

Note that Lemmas \ref{tri id lemma}, \ref{rect id lemma} state an equality of
rational functions.
Therefore, the equality still holds if we choose to
expand the RHS into power series in any consistent way.
There are four choices for the expansion of \Ref{rect id},
depending on whether $x,y$ are taken
in the neighborhood of $0$ or $\infty$.
For \Ref{tri id} there are 6 ways to expand
since
there is an additional choice of $xy$
being in the neighborhood of $0$ or $\infty$.

The lemmas \ref{tri id lemma}, \ref{rect id lemma} are simple examples
of writing a character of a convex polytope
as a sum over vertices of characters of the corresponding cones.
The second lemma in fact follows from an even more simple example,
namely from the expression for a 1-dimensional segment:
\be
\sum_{a\leq m\leq b} x^m=\frac{x^a}{1-x}+\frac{x^b}{1-x^{-1}}.
\ee
We remark that the Weyl formula for the characters of
irreducible finite dimensional modules over semisimple
Lie algebras have a similar structure.


\section{The extremal operators}
Let $P,Q,R$  be monomials of the form
$q^{\al}z_1^{\beta}z_2^{\gamma}$ with $\al,\beta,\gamma\in\Z_{\geq 0}$. 
We denote by $[P,Q,R]$ the $(k+1)(k+2)/2$-dimensional vector whose 
$(i,l)$-th component is given by 
\be
[P,Q,R]_{i,l}=[P(1:{Q\over P}:{R \over P})]_{i,l}
=P^{k-l}Q^{l-i}R^i.
\ee
In this notation we have 
\bean\label{init}
\chi_{i,l}^{(0)}(q,z_1,z_2)&=&\delta_{i,0}\delta_{l,0}=[1,0,0]_{i,l},\notag\\
\chi_{i,l}^{(1)}(q,z_1,z_2)&=&z_2^i \delta_{l,i}=[1,0,z_2]_{i,l},\\
\chi_{i,l}^{(2)}(q,z_1,z_2)&=&\sum_{j=0}^{k-l}q^{l-i+j}z_2^{l+j}=
\frac{[1,qz_2,z_2]_{i,l}}{1-qz_2}+\frac{[qz_2,z_2,z_2]_{i,l}}{1-(qz_2)^{-1}}. 
\notag
\eean

Consider the linear span $V$ of vectors 
\begin{eqnarray}
f[P,Q,R],  \qquad f\in \C[[q,z_2]]((z_1^{-1})). 
\label{simple}
\end{eqnarray}
Normally we will write an element of $V$ simply as $v$,   
without exhibiting the dependence on $q,z_1,z_2$ explicitly. 
We say $v$ is {\em simple} if it has the form \eqref{simple}. 
We call $f$ and $[P,Q,R]$ 
the {\em scalar part} and the {\em vector part}, respectively. 
Since $f[aP,aQ,aR]=fa^k[P,Q,R]$ holds for a monomial $a$,  
the scalar part 
and the vector part are determined up to this freedom. 

For a vector-valued function $g(q,z_1,z_2)$, let 
\bean\label{shift}
(Sg)(q,z_1,z_2)=g(q,z_1,qz_2)
\eean
denote the $q$-shift operator in $z_2$. 
We introduce the {\em extremal operators} $A,B,C,D,E$ by the
formulas 
\bea
A(f[P,Q,R])&=&
S\left(\frac{f[P,Q,q^{-1}z_2P]}{(1-z_1z_2Q/P)(1-R/Q)}\right),\\
B(f[P,Q,R])&=&
S\left(\frac{f(1-Q/P)[P,R,q^{-1}z_2P]}{(1-z_1z_2Q/P)(1-Q/R)(1-R/P)}\right),\\
C(f[P,Q,R])&=&
S\left(\frac{f[z_1z_2Q,Q,q^{-1}z_2P]}
{(1-(z_1z_2)^{-1}P/Q)(1-R/Q)}\right),\\
D(f[P,Q,R])&=&
S\left(\frac{(1-(z_1z_2)^{-1})f[z_1z_2Q,R,q^{-1}z_2P]}
{(1-(z_1z_2)^{-1}P/Q)(1-(z_1z_2)^{-1}R/Q)(1-Q/R)}\right),\\
E(f[P,Q,R])&=&
S\left(\frac{f[R,R,q^{-1}z_2P]}{(1-z_1z_2Q/R)(1-P/R)}\right),
\eea
extended by linearity. 
The meaning of the RHS is as follows. 
For a monomial $X=q^\alpha z_1^\beta z_2^\gamma$ 
($\alpha,\gamma\in\Z_{\ge 0}, \beta\in\Z$) 
such that $\alpha+\gamma>0$ or $\alpha=\gamma=0,\beta<0$, 
$1/(1-X^{\pm 1})$ means 
its expansion in non-negative powers of $X$.  
We define an operator $G(=A,B,C,D,E)$ on $f[P,Q,R]$ 
when the denominators appearing in the RHS are all 
of this form. 
Otherwise we do not define $G(f[P,Q,R])$.  

The main point of introducing the operators $A,B,C,D,E$ is:
\begin{lemma}
Suppose operators $A,B,C,D,E$ are defined on a simple vector $v$,  
and let $M(q,z_1,z_2)$ be the matrix \Ref{matrix}.
Then we have 
\be
M(q,z_1,z_2)v(q,z_1,qz_2)=((A+B+C+D+E)v)(q,z_1,z_2).
\ee
\end{lemma}
\begin{proof}
It suffices to consider the case $v=[P,Q,R]$. 
Written out explicitly, 
the $(i,l)$-component of the LHS reads
\begin{eqnarray*}
&&\sum_{l-i\le i'\le l'\le k-i}(qz_1z_2)^{l'-i'}z_2^i \,
S(P)^{k-l'}S(Q)^{l'-i'}S(R)^{i'}
\\
&&\quad+
\sum_{i'< l-i\le l'\le k-i}(qz_1z_2)^{l'-l+i}z_2^i \,
S(P)^{k-l'}S(Q)^{l'-i'}S(R)^{i'}.
\end{eqnarray*}
We apply Lemmas \ref{tri id} and \ref{rect id}
to obtain seven terms corresponding to vertices in 
Figure \ref{pentagon}.   
Namely 
\bean\label{inter}
M(q,z_1,z_2)v(q,z_1,qz_2)=
(\left(B_1+E+D_1+A+B_2+D_2+C\right)v)(q,z_1,z_2),
\eean
where
\bea
B_1([P,Q,R])&=&
S\left(\frac{[P,R,q^{-1}z_2P]}{(1-z_1z_2Q/P)(1-R/P)}\right),\\
B_2([P,Q,R])&=&
S\left(\frac{-[P,R,q^{-1}z_2P]}{(1-z_1z_2Q/P)(1-R/Q)}\right),\\
D_1([P,Q,R])&=&
S\left(\frac{[z_1z_2Q,R,q^{-1}z_2P]}
{(1-(z_1z_2)^{-1}P/Q)(1-(z_1z_2)^{-1}R/Q)}\right),\\
D_2([P,Q,R])&=&
S\left(\frac{-[z_1z_2Q,R,q^{-1}z_2P]}
{(1-(z_1z_2)^{-1}P/Q)(1-R/Q)}\right).
\eea
The lemma follows from the identities $B=B_1+B_2$, $D=D_1+D_2$. 
\end{proof}

If the operators $A,B,C,D,E$ were always defined, 
the lemma could be applied repeatedly to write the character as 
\bean\label{trivial}
\chi^{(N)}(q,z_1,z_2)&=&
(A+B+C+D+E)^{N-1}[1,0,z_2],\notag\\
\chi(q,z_1,z_2)&=&
(A+B+C+D+E)^\infty[1,0,z_2].
\eean
Unfortunately it is not the case.  
For example,
\be
CBCAE[P,Q,R]=f(P,Q,R)[q^3z_1z_2^2\tilde R,q^2z_2\tilde
R,q^3z_1z_2^2\tilde R], \qquad \tilde R=S^5(R).
\ee
Therefore the operators $BCBCAE$ and $ECBCAE$ 
are not defined because $0$ is produced in the denominator. 
(Note however that the sum $(B+E)CBCAE$ is well defined.)
This example shows that 
there is no subspace $W$ which contains $[1,0,0]$ and
is stable under $A,B,C,D,E$.
Such a difficulty does not arise 
in the case treated in \cite{FJLMM}.

Nevertheless it is possible to express 
$\chi(q,z_1,z_2)$ as a sum over a subset of non-commutative 
monomials in $A,B,C,D,E$, 
with all terms being defined in the sense above. 
In Theorem \ref{main} below we write down the formula and 
show directly that it is the unique 
solution of the difference equation \eqref{diff eqn}.
Informally speaking, it means that 
the rest of the terms in \eqref{trivial} 
cancel with each other, including the non-defined ones.  
In the next two sections we 
analyze the mechanism of the cancellation and  
explain how the formula was found. 
This part is meant to motivate Theorem \ref{main},   
although it is logically unnecessary.

\section{Summation graph}
In this section we prepare some language for 
Section \ref{inf cancel}.

A {\em monomial} $\mc M$ is an ordered composition of operators
$G_1G_2\dots G_n$, where $G_i\in\{A,B,C,D,E\}$.
We call $n$ the degree of $\mc M$.

Now, our goal is to find different monomials which, 
when acted upon a simple vector $v\in V$, give rise to
the same vector part.
For example, the vector parts of $BC([P,Q,R])$ and 
$BD([P,Q,R])$ are
both equal to $S^2([z_1z_2Q,q^{-1}z_2P,z_1z_2^2q^{-2}Q])$. 
Therefore 
the sum $B(C+D)$ applied to a simple vector is again a simple
vector. 
Note that even though some monomials may not be defined,
the action of a monomial on the vector part is always 
defined. 

Introduce the maps $\sigma_G:\;\{1,2,3\}\to\{1,2,3\}$, where
$G\in\{A,B,C,D,E\}$, by the formulas
\bea
\ol\sigma_A=(1,2,1),\qquad \ol\sigma_B=(1,3,1), 
\qquad \ol\sigma_C=(2,2,1),\\
\ol\sigma_D=(2,3,1),\qquad \ol\sigma_E=(3,3,1),
\eea
where $\ol\sigma_G=(\sigma_G(1),\sigma_G(2),\sigma_G(3))$.

We define the map $\sigma_{\mc M}$ for all monomials ${\mc M}$
by the product rule
\be
\sigma_{{\mc M}_1{\mc M}_2}=\sigma_{{\mc M}_2}\sigma_{{\mc M}_1}.
\ee
In other words $\sigma$ is an anti-homomorphism of the semi-group of
non-commutative monomials in $A$,$B$,$C$,$D$,$E$
to the semi-group of endomorphisms of $\{1,2,3\}$.

The maps $\sigma_{\mc M}$ describe how $\mc M$ permutes $P,Q,R$ in the vector
part. Namely,
\be
{\mc M}([P_1,P_2,P_3])=S^{\deg \mc
M}(f[u_1 P_{\sigma(1)},u_2 P_{\sigma(2)},u_3 P_{\sigma(3)}]),
\ee
where $u_1,u_2,u_3$ denote factors independent of $P_1,P_2,P_3$.

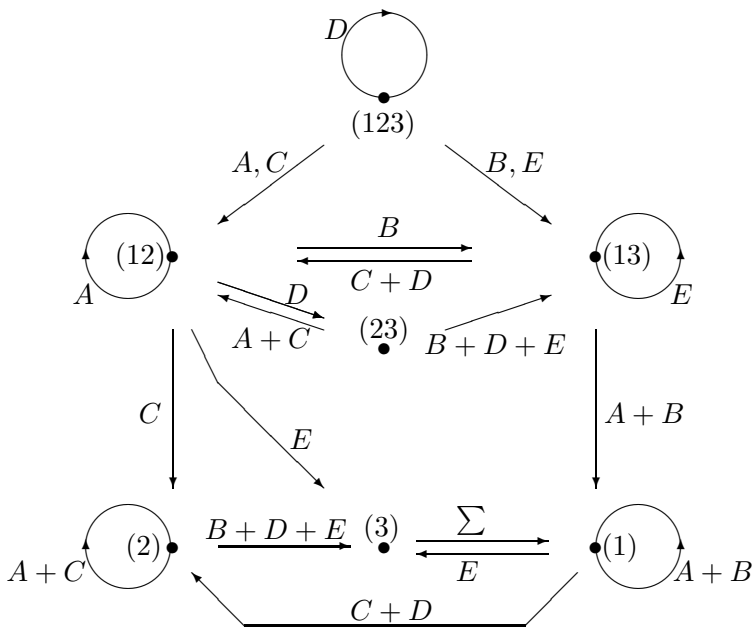
\begin{figure}
\begin{picture}(0,240)(0,0)
\put(0,200){$\bullet$}
\put(-80,140){$\bullet$}
\put(80,140){$\bullet$}
\put(0,105){$\bullet$}

\put(-80,30){$\bullet$}
\put(80,30){$\bullet$}
\put(0,30){$\bullet$}

\put(3,219){\circle{30}}
\put(6,235){\vector(1,0){0}}

\put(-94,143){\circle{30}}
\put(-110,146){\vector(0,1){0}}
\put(99,143){\circle{30}}
\put(115,146){\vector(0,1){0}}

\put(-94,33){\circle{30}}
\put(-110,36){\vector(0,1){0}}
\put(99,33){\circle{30}}
\put(115,36){\vector(0,1){0}}

\put(-10,190){$(123)$}
\put(-99,140){$(12)$}
\put(85,140){$(13)$}
\put(-7,112){$(23)$}
\put(85,30){$(1)$}
\put(-5, 37){$(3)$}
\put(-95,30){$(2)$}

\put(-20,185){\vector(-4,-3){40}}
\put( 26,185){\vector( 4,-3){40}}

\put( 36,141){\vector(-1,0){66}}
\put(-30,146){\vector( 1,0){66}}

\put(-60,133){\vector(3,-1){40}}
\put(-20,115){\vector(-3,1){40}}
\put( 26,115){\vector(3,1){40}}
\put(-60,133){\vector(3,-1){40}}

\put(-70,115){\line(1,-2){10}}
\put(-60,95){\vector(1,-1){40}}
\put(-77,115){\vector(0,-1){60}}
\put( 83,115){\vector(0,-1){60}}

\put(-60,33){\vector(1,0){50}}
\put( 15,35){\vector(1,0){50}}
\put( 65,30){\vector(-1,0){50}}

\put(-50, 03){\vector(-1,1){20}}
\put( 56, 03){\line(1,1){20}}
\put( -50, 03){\line(1,0){106}}

\put(-20,225){$D$}
\put(-55,175){$A,C$}
\put( 41,175){$B,E$}

\put(0,150){$B$}
\put(-10,130){$C+D$}
\put(-55,108){$A+C$}
\put(-35,125){$D$}
\put(18,106){$B+D+E$}

\put(86,80){$A+B$}
\put(-90,80){$C$}
\put(-33,70){$E$}

\put(-65,36){$B+D+E$}
\put(30,40){$\sum$}
\put(30,20){$E$}
\put(-10,5){$C+D$}

\put(-115,125){$A$}
\put(111,125){$E$}
\put(-140,20){$A+C$}
\put(112,20){$A+B$}

\end{picture}
\caption{Summation graph (here $\sum=A+B+C+D+E$).}\label{add graph}
\end{figure}

Introduce a directed graph $\Gamma$, 
which we call {\em summation graph}
(see Figure \ref{add graph}).
The vertices 
of $\Gamma$ are non-empty subsets $I$ of $\{1,2,3\}$. Each
vertex emits 5 arrows, labeled $A,B,C,D,E$.
The arrow with label $G$ which starts from a vertex $I$
ends at the vertex $\sigma_G(I)$.

Note that the summation graph $\Gamma$ has
three {\em floors} consisting of subsets $I$ of given
cardinality. We picture vertices corresponding to subsets of larger
cardinality above those correspodning to subsets of smaller
cardinality.  Then no arrow goes ``up''.

In our picture, arrows with the same beginning and end 
are represented by one arrow with several labels.  
For arrows $G_1,G_2$ with a common source $I\subset\{1,2,3\}$, 
we write $G_1+G_2$ if and only if 
the $i$-th component of the 
vector parts of $G_1([P,Q,R])$ and $G_2([P,Q,R])$ 
coincide for all $i\in I$. 
For example, the vector parts of $A([P,Q,R])$ and $C([P,Q,R])$ 
coincide in the second and the third components. 
Therefore we write $A+C$ above the arrows coming from 
$(2)$ and $(2,3)$, but do not do so for $(1,2,3)$.

We identify monomials in $A,B,C,D,E$ with (oriented) paths in $\Gamma$
starting at $(1,2,3)$, by reading a monomial $\mc M$ from left to right
and choosing the corresponding arrows
accordingly. For example, the monomial $ACB$
corresponds to the path
\be
(1,2,3){\buildrel A\over\longrightarrow} (1,2)
{\buildrel C\over\longrightarrow}(2){\buildrel B\over\longrightarrow}(3).
\ee

The vector part of $\mc M([P_1,P_2,P_3])$ corresponding to a path
ending at a vertex $I$ depends only on $P_i$ with $i\in I$.
Therefore we have the following obvious lemma.
\begin{lemma}
Let $\mc M$ be a monomial.
Suppose that the corresponding path in the summation graph
ends at a source of an arrow $G_1+\dots+G_s$, where $G_i$ are
distinct elements of $\{A,B,C,D,E\}$.
Then the operator $\mc M(G_1+\dots+G_s)$ maps 
simple vectors $v\in V$ to simple vectors
if it is defined. 
\end{lemma}

We will use such summations in Section \ref{inf cancel}
to combine several terms in \Ref{trivial} in one.
It turns out that the sum of the corresponding scalar parts
in our computations is completely factorized to
``linear'' factors.
We warn the reader that 
this is not necessarily the case
in general (e.g. in the bottom floor for finite $N$).  

\section{Cancellation in the case $N=\infty$}\label{inf cancel}
{}From now on we concentrate on the case $N=\infty$.
The purpose of this section is to describe
the structure of the resulting formula and the way
how it was obtained.
In the argument below, we assume that 
all compositions of operators are defined.
After finding the final formula 
we turn to proving it by direct means.

In formula \Ref{trivial}, the
RHS is a sum of monomials in $A,B,C,D,E$ of infinite
degree applied to $[1,0,z_2]$.
Our first step is to choose another 
vector in place of $[1,0,z_2]$, 
so that it suffices 
to deal only with monomials of finite degree. 

Let $\mc M$ be an arbitrary monomial of infinite degree
and $\mc N$ of finite degree.
We observe that $\mc M G \mc N ([1,0,z_2])=0$ holds for $G=C,D,E$.

The reason is as follows.
We claim, that if $ G \mc N ([1,0,z_2])\neq 0$, then all
components of its 
vector part contain a factor $z_2$. 
Indeed, if $\mc N=1$ then
$C([1,0,z_2])=D([1,0,z_2])=0$, 
and the vector part of $E([1,0,z_2])$
is $[qz_2,qz_2,z_2]$. If $\mc N\neq 1$, then
the only component which does not 
depend on $z_2$ in ${\mc N}([1,0,z_2])$ 
may be the first one, and our claim follows.
Then $[z_2P,z_2Q,z_2R]=z_2^k[P,Q,R]$ 
has an overall factor $z_2^k$.
It vanishes when $\mc M$ is applied, 
because of an infinite shift $S^\infty$.

The formula \Ref{trivial} would then take the form
\bean \label{finite paths}
\chi(q,z_1,z_2)=
\sum_{\deg{\mc M}<\infty}\mc M(A+B)^\infty([1,0,z_2]).
\eean
In the language of paths 
in the summation graph $\Gamma$, 
the formula \Ref{finite paths} means that
we sum over all paths of finite length which
start at $(1,2,3)$ and end at $(1)$.

We compute $(A+B)^\infty([1,0,z_2])$ using the identity
\be
(A+B)^n=A^n+\sum_{j=0}^{n-1}A^{j}B(A+B)^{n-j-1}.
\ee
We have explicitly
\be
A^nB(A+B)^m[(1,0,z_2)]=\frac{[1,q^{n+1}z_2,z_2]}
{(1-q^{n+1}z_2)\prod_{j=-n\atop j\not=0}^{m-1}(1-q^j)
\prod_{j=n+2\atop j\not=2n+2}^{2n+m+1}(1-q^jz_1z_2^2)},
\ee
and $A([1,0,z_2])=0$. Taking the limit $n\to\infty$ we obtain the
following definition.

Introduce the vector $v_\infty\in V$ by the formulas
\begin{eqnarray}\label{vinf1}
v_\infty=\sum_{n=1}^\infty f_nv_n,
\end{eqnarray}
where
\begin{eqnarray}\label{vinf2}
v_n=[1,q^nz_2,z_2], \qquad
f_n=\frac{(-1)^{n-1}q^{n(n-1)/2}(1-q^{2n}z_1z_2^2)}
{(q)_\infty(q^{n+1}z_1z_2^2)_\infty(q)_{n-1}(1-q^nz_2)}.
\end{eqnarray}
We have an identity $(A+B)v_\infty=v_\infty$ (see Lemma \ref{stable}
below).

Now, we describe the mechanism of the cancellation and the 
structure of the final formula. 

First of all we infer, on the basis of examples, 
that paths 
which pass through vetrices $(2)$ or $(3)$ cancel out. 
In other words the surviving monomials correspond to paths in the
summation graph which never go to the bottom level until the tail
$(A+B)^\infty$ is reached.
We have no direct proof of this statement. 

After that reduction, we still have several 
cycles which paths can wrap around in the middle floor. 
Let us set 
\be
L=(C+D)D(B+D+E), 
\quad
\bar{L}=D(B+D+E)(C+D).  
\ee
We have several operator identities which explain a part of the 
cancellation. 
For example, we have 
\be
BE=0,
\ee
which implies in particular that 
if $E$ cycle appears after $L$ then the
corresponding term is zero. 
We have also the identities 
valid for any $m\in\Z_{\geq 0}$,
\begin{eqnarray*}
&&X \bar{L}^m D(A+C)=0, \quad X\bar{L}^m B=0
\qquad (X=A,C), 
\\
&&YL^m(C+D)D(A+C)=0,
 \quad Y{L}^m(C+D)B=0 \qquad (Y=B,E),  
\end{eqnarray*}
which imply that the cycle $D(A+C)$,
 $(C+D)B$ give zero contributions.
(We stress that we do not use any of these operator 
identities for our proofs.)
So, we are left with the cycles $A$ (at the vertex $(1,2)$),
$E$ (at the vertex $(1,3)$)
and $L$ (starting at the vertex $(1,3)$). 

Finally, there is an ordering of the cycles 
in the following sense.
We call a path in the 
summation graph $\Gamma$ {\em good} 
if the following 3 conditions are fulfilled:
\begin{itemize}
\item it ends at the vertex $(1,3)$, 
\item it does not contain both $A$ and $E$ cycles,
\item it enters neither $A$ nor $E$ cycles after passing through
$L$ cycle.
\end{itemize}
A path which is not good is called {\em bad}.
We infer that in the formula \Ref{finite paths}
bad paths cancel out and only good paths survive.

These considerations lead us to the following theorem.
\begin{theorem}\label{main}
We have the following identity
\bean\label{main eq}
\chi(q,z_1,z_2)=\sum_{n,m,s=0}^\infty
\left(D^nA^{m}BL^sv_\infty
+ D^nCA^{m}BL^sv_\infty +
D^nE^{m+1}L^sv_\infty +\right.\notag\\
\left. +D^nCA^{m}D(B+D+E)L^sv_\infty+ D^nA^{m+1}D(B+D+E)L^sv_\infty\right),
\eean
where $v_\infty$ is defined in \eqref{vinf1},\eqref{vinf2}.
All terms in the RHS are defined. 
\end{theorem}

\begin{proof}
For the proof we use the explicit 
formulas for each term in the RHS
given in Section \ref{expl section}. 
By a direct computation we verify that they are all defined, 
and that each operator $A,B,C,D,E$ are also defined on them. 

Substitution $z_2=0$ to the RHS clearly gives $0$ (because of the
nontrivial vector part) unless $l=i=0$, when we get $1$ 
(coming from the term $Bv_\infty$).
By the uniqueness part of Lemma \ref{existence},  
it is enough to show 
that multiplication by $A+B+C+D+E$ does not change the RHS.

There are two things to show. 
First, we have to show that  
multiplication by $A+B+C+D+E$ 
reproduces all the terms, and 
second, that the extra terms appearing cancel out.

The first statement, informally speaking, 
can be reformulated as follows.
If a monomial $G\mc M$ with $G\in\{A,B,C,D,E\}$
is good, then monomial $\mc M$ itself is also good. 
We verify it by case checking. The nontrivial 
cases are $BL^sv_\infty$, $EL^sv_\infty$, $CD^2L^sv_\infty$ and
$AD^2L^sv_\infty$ with $s\in\Z_{>0}$. We have to show that
$L^sv_\infty$ and $D^2L^sv_\infty$ occur in the RHS.
Let us consider the case of $L^sv_\infty$, the other case is similar.

Expanding $L$ into a sum of monomials, 
we see that all monomials in
$L^sv_\infty$ are good, except for the monomial
$D^{3s}v_\infty$, which we replace by the sum of good
monomials $\sum_{n=0}^\infty D^{3s}A^nBv_\infty$, using the identity
(see Lemma \ref{stable})
\be
(1-A)^{-1}Bv_\infty=\sum_{n=0}^\infty A^nBv_\infty=v_\infty.
\ee

For the second statement, we claim that cancellations always happen in
pairs. We list the necessary equalities which are checked using the
explicit formulas given in Section \ref{expl section}.

In all equalities below, $m,n,s$ stand for arbitrary non-negative integers.
First we have cancellations of bad monomials of type $A\mc M$, where
$\mc M$ is a good monomial.
\bea
(AD^{3n+3}A^mBL^s+AD^{3n+2}CA^mBL^s)v_\infty=0,\\
(AD^{3n+1}A^{m+1}BL^s+AD^{3n+1}CA^mBL^s)v_\infty=0,\\
(AD^{3n+2}A^{m}BL^{s+1}+AD^{3n+2}A^{m+1}D(B+D+E)L^s)v_\infty=0,
\eea
\bea
(AD^{3n}CA^mBL^s+AD^{3n}CA^mD(B+D+E)L^s)v_\infty=0,\\
(AD^{3n}E^{2m+1}L^s+AD^{3n}E^{2m+2}L^s)v_\infty=0, \\
(AD^{3n+1}E^{2m+3}L^s+AD^{3n+2}E^{2m+2}L^s)v_\infty=0,\\
\eea
\bea
(AD^{3n+2}E^{2m+1}L^s+AD^{3n+1}E^{2m+2}L^s)v_\infty=0,\\
(AD^{3n+3}A^{m+1}D(B+D+E)L^s+AD^{3n+2}CA^{m+1}D(B+D+E)L^s)v_\infty=0,\\
(AD^{3n+1}A^{m+1}D(B+D+E)L^s+AD^{3n+1}CA^{m}D(B+D+E)L^s)v_\infty=0.
\eea
Then another set of $9$ equalities are obtained from the above $9$
by replacing the first factor $A$ by $C$.
They describe cancellations of bad monomials of type $C\mc M$, where
$\mc M$ is a good monomial.

In the case $B\mc M$, where $\mc M$ is a good monomial, we have the
following equalities.
\bea
(BD^{3n}A^mBL^{s+1}+BD^{3n}A^{m+1}D(B+D+E)L^s)v_\infty=0,\\
(BD^{3n+1}A^{m}BL^s+BD^{3n}CA^mBL^s)v_\infty=0,\\
(BD^{3n+2}A^{m+1}BL^{s}+BD^{3n+2}CA^mBL^s)v_\infty=0,
\eea
\bea
(BD^{3n+1}CA^mBL^s+BD^{3n+1}CA^{m}D(B+D+E)L^s)v_\infty=0,\\
(BD^{3n+3}E^{2m+1}L^s+BD^{3n+2}E^{2m+2}L^s)v_\infty=0, \\
(BD^{3n+1}E^{2m+1}L^s+BD^{3n+1}E^{2m+2}L^s)v_\infty=0,\\
\eea
\bea
(BD^{3n+2}E^{2m+3}L^s+BD^{3n+3}E^{2m+2}L^s)v_\infty=0,\\
(BD^{3n+1}A^{m+1}D(B+D+E)L^s+BD^{3n}CA^{m+1}D(B+D+E)L^s)v_\infty=0,\\
(BD^{3n+2}A^{m+1}D(B+D+E)L^s+BD^{3n+2}CA^{m}D(B+D+E)L^s)v_\infty=0.
\eea
The cancellations of bad monomials of type $E\mc M$, where
$\mc M$ is a good monomial are given by replacing first factor $B$
in the above formulas by $E$.

Note that if $\mc M$ is good then $D\mc M$ is also good.
\end{proof}
Note that each term in the RHS of Theorem \ref{main}
is a product of a monomial in $q,z_1,z_2$  with ``linear'' factors of
the form $(1-q^{s+1})^{\pm 1}$, $(1-q^{s}z_1^{-1})^{\pm 1}$,
$(1-q^{s}z_2)^{\pm 1}$,
$(1-q^{s}z_1z_2)^{\pm 1}$, $(1-q^{s}z_1z_2^2)^{\pm 1}$ with $s\geq 0$.

Also note, that in order to prove the cancellation (or in order to
write closed formulas) one has to break 5 families in the
RHS of \Ref{main eq} into 18 families.
Then each function on the RHS is a power series in
$q$ whose coefficients are
formal Laurent series in 
$z_1^{-1}$ over the ring of polynomials $\Z[z_2]$.
We know by corollary \ref{integercoeff} that 
terms containing $z_1^{-1}$ disappear after the summation. 
	
Finally, we remark that some cancellations of the above happen between
monomials of different degrees, for example
$CEv_\infty+CE^2v_\infty=0$. However, we could rewrite this equality
as a sum of monomials of the same degree:
$CE(A+B)v_\infty+CE^2v_\infty=0$.

\section{The explicit formula}\label{expl section}

We start with an identity which plays a crucial role in
our computations:
\bean\label{id}
\frac{(q)_\infty}{(qx)_\infty (qy)_\infty}=\sum_{n=1}^\infty\frac{(-1)^{n-1}q^{n(n-1)/2}(1-q^{2n}xy)}{(q)_{n-1}(q^{n+1}xy)_\infty(1-q^nx)(1-q^ny)}.
\eean
In fact this is a special case
of Jackson's ${}_6\Phi_5$ formula which states
(see for example \cite{SlBk} p.102, formula (3.4.2.3)):
\be
\sum_{n=0}^\infty
\frac{(1-q^{2n}a)(a)_n(b)_n(c)_n(d)_n(qa/bcd)^n}{(1-a)(q)_n(qa/b)_n(qa/c)_n(qa/d)_n}
=
\frac{(qa)_\infty(qa/cd)_\infty(qa/bd)_\infty(qa/bc)_\infty}
{(qa/b)_\infty(qa/c)_\infty(qa/d)_\infty(qa/bcd)_\infty}.
\ee
Identity \Ref{id} is obtained by setting $a=xy$, $b=x$, $c=y$, $d=\infty$.

Here is an example.
\begin{lemma}\label{stable}
We have $(A+B)v_\infty=v_\infty$. 
Moreover, $A^nf_1v_1=f_{n+1}v_{n+1}$
and $Bv_\infty=f_1v_1$.
\end{lemma}
\begin{proof} The identity $A^nf_1v_1=f_{n+1}v_{n+1}$ is obvious. The
identity $Bv_\infty=f_1v_1$ is equivalent to \Ref{id} with
$x=q^2z_1z_2^2$ and $y=1$.
\end{proof}

Note that $v_\infty$ is a sum of infinitely many simple
vectors with vector parts $[1,q^nz_2,z_2]$ and the first
component of these vector parts is always $1$. Therefore, if $\mc M$
is a monomial corresponding to a path which ends at
vertex $(1)$ in the summation graph, then $\mc M[1,q^nz_2,z_2]$ does not
depend on $n$, and $\mc Mv_\infty$ is a simple vector.
In particular, all formulas below are simple vectors. The
summation of scalar parts can be always explicitly performed using \Ref{id}.

Now we finish with explicit formulas for the terms in the RHS of the formula
in Theorem \ref{main}. These formulas are obtained by 
an explicit computation, 
induction on $m,n,s$ and using formula \Ref{id}.

We define quadratic forms $\alpha$, $\beta$, $\gamma$ and $\delta$  by
\bea
&&\alpha_{n,m,s}=3(n+s)^2+2ms,\\
&&\beta_{n,m,s}=3(n+s)^2+m^2+4ms+3mn,\\
&&\gamma_{n,m,s}
=\frac{21}{2}n^2+\frac{1}{2}m^2+\frac{13}{2}s^2+3nm+15ns+5ms,
\\
&&\delta_{n,m,s}
=\frac{11}{2}n^2+\frac{3}{2}m^2+5 s^2+5nm+10ns+6ms.
\eea
Then we have the following formulas.

\begin{align*}
&D^{3n-2}A^mBL^sv_\infty=(-1)^{n+m+s}(z_1z_2^3)^{n+s}z_2^{-1}
q^{\gamma_{n,m,s}-\frac{17}{2}n-\frac{1}{2}m-\frac{7}{2}s+2}\times
\\
&\times
\frac{(qz_1z_2)_{2n-1}(q^{3n+m+s}z_1z_2)_s(1-q^{6n+2m+2s-2}z_1z_2^2)}
{(q)_\infty(q)_{2n-2}(q)_{2n+m+2s-1}(q^{3n+m+s-1}z_2)_{n+s}(q^{2n+m+2s}z_1z_2)_n(q^{4n+m+2s-1}z_1z_2^2)_\infty}\times
\\
&\times
[q^{\al_{n,m,s}+m}(z_1z_2^2)^{n+s}
(1:q^{-2n-m-2s}(z_1z_2)^{-1}:q^{1-4n-m-2s}(z_1z_2)^{-1})],\;
n\ge 1;m,s\ge 0,
\end{align*}

\begin{align*}
&D^{3n-1}A^mBL^sv_\infty=
(-1)^{n+m+s}(z_1z_2^3)^{n+s}z_2^{-1}
q^{\gamma_{n,m,s}-\frac{7}{2}n-\frac{1}{2}m-\frac{1}{2}s}\times\\
&
\times
\frac{(qz_1z_2)_{2n-1}(q^{3n+m+s+1}z_1z_2)_s(1-q^{6n+2m+2s}z_1z_2^2)}
{(q)_\infty(q)_{2n-1}(q)_{2n+m+2s}(q^{3n+m+s}z_2)_{n+s}
(q^{2n+m+2s+1}z_1z_2)_n(q^{4n+m+2s+1}z_1z_2^2)_\infty}\times
\\
&\times
[q^{\al_{n,m,s}}(z_1z_2^2)^{n+s}
(1:q^{-2n}(z_1z_2)^{-1}:q^{2n+m+2s}z_2)], \qquad \quad n\ge 1,m\ge 0,s\ge 0,
\end{align*}

\begin{align*}
&D^{3n}A^{m-1}BL^sv_\infty=(-1)^{n+m+s+1}(z_1z_2^3)^{n+s}
q^{\gamma_{n,m,s}+\frac{1}{2}n-\frac{1}{2}m-\frac{1}{2}s}\times
\\
&\times
\frac{(qz_1z_2)_{2n}(q^{3n+m+s+1}z_1z_2)_s(1-q^{6n+2m+2s}z_1z_2^2)}
{(q)_\infty(q)_{2n}(q)_{2n+m+2s-1}(q^{3n+m+s}z_2)_{n+s+1}
(q^{2n+m+2s+1}z_1z_2)_n(q^{4n+m+2s+1}z_1z_2^2)_\infty}\times
\\
&\times
[q^{\alpha_{n,m,s}}(z_1z_2^2)^{n+s}
(1:q^{4n+m+2s}z_2:q^{2n}z_2)], \qquad \qquad \qquad n\ge 0,m\ge 1,s\ge 0,
\end{align*}

\begin{align*}
&D^{3n-2}CA^{m-1}BL^sv_\infty=
(-1)^{n+m+s}(z_1z_2^3)^{n+s}z_2^{-1}
q^{\gamma_{n,m,s}-\frac{13}{2}n-\frac{1}{2}m-\frac{7}{2}s+1}\times
\\
&\times
\frac{(qz_1z_2)_{2n-1}(q^{3n+m+s}z_1z_2)_s(1-q^{6n+2m+2s-2}z_1z_2^2)}
{(q)_\infty(q)_{2n-1}(q)_{2n+m+2s-1}(q^{3n+m+s-1}z_2)_{n+s}
(q^{2n+m+2s}z_1z_2)_n(q^{4n+m+2s}z_1z_2^2)_\infty}\times
\\
&\times
[q^{\al_{n,m,s}+m}(z_1z_2^2)^{n+s}
(1:q^{-2n-m-2s}(z_1z_2)^{-1}:q^{2n+m+2s}z_2)],\;\; n\ge 1,m\ge 1,s\ge 0,
\end{align*}

\begin{align*}
&D^{3n-1}CA^{m-1}BL^sv_\infty=
(-1)^{n+m+s}(z_1z_2^3)^{n+s}
q^{\gamma_{n,m,s}+\frac{1}{2}n+\frac{1}{2}m+\frac{3}{2}s}\times
\\
&\times
\frac{(qz_1z_2)_{2n}(q^{3n+m+s+1}z_1z_2)_s(1-q^{6n+2m+2s}z_1z_2^2)}
{(q)_\infty(q)_{2n-1}(q)_{2n+m+2s}(q^{3n+m+s}z_2)_{n+s+1}(q^{2n+m+2s+1}z_1z_2)_n(q^{4n+m+2s+1}z_1z_2^2)_\infty}
\times
\\
&\times
[q^{\alpha_{n,m,s}}(z_1z_2^2)^{n+s}
(1:q^{4n+m+2s}z_2:q^{2n+m+2s}z_2)], \qquad \qquad n\ge 1,m\ge 1,s\ge 0,
\end{align*}

\begin{align*}
&D^{3n-3}CA^{m}BL^sv_\infty=
(-1)^{n+m+s+1}(z_1z_2^3)^{n+s}z_2^{-1}
q^{\gamma_{n,m,s}-\frac{17}{2}n-\frac{1}{2}m-\frac{7}{2}s+2}
(q)_\infty^{-1}\times
\\
&\times
\frac{(qz_1z_2)_{2n-2}(q^{3n+m+s}z_1z_2)_s(1-q^{6n+2m+2s-2}z_1z_2^2)}
{ (q)_{2n-2}(q)_{2n+m+2s-1}(q^{3n+m+s-1}z_2)_{n+s}
(q^{2n+m+2s+1}z_1z_2)_{n-1}(q^{4n+m+2s-1}z_1z_2^2)_\infty}\times
\\
&\times
[q^{\al_{m,n,s}+m}(z_1z_2^2)^{n+s}
(1:q^{-2n+1}(z_1z_2)^{-1}:q^{2-4n-m-2s}(z_1z_2)^{-1})],\;n\ge 1,m\ge 0,s\ge 0,
\end{align*}

\begin{align*}
&D^{3n-1}A^{m+1}D(B+D+E)L^{s-1}v_\infty=(-1)^{n+m+s}(z_1z_2^3)^{n+s}z_2^{-1}
q^{\gamma_{n,m,s}-\frac{11}{2}n-\frac{3}{2}m-\frac{5}{2}s}\times
\\
&\times
\frac{(qz_1z_2)_{2n-1}(q^{3n+m+s+1}z_1z_2)_{s-1}(1-q^{6n+2m+2s}z_1z_2^2)}
{(q)_\infty(q)_{2n-1}(q)_{2n+m+2s-1}(q^{3n+m+s}z_2)_{n+s}
(q^{2n+m+2s+1}z_1z_2)_{n-1}(q^{4n+m+2s}z_1z_2^2)_\infty}\times
\\
&\times
[q^{\al_{n,m,s}}(z_1z_2^2)^{n+s}
(1:q^{-2n}(z_1z_2)^{-1}:q^{-4n-m-2s}(z_1z_2)^{-1})],\qquad
\quad n\ge 1,m\ge 0,s\ge 1,
\end{align*}

\begin{align*}
&D^{3n-2}A^{m}D(B+D+E)L^{s}v_\infty=
(-1)^{n+m+s+1}(z_1z_2^3)^{n+s}
q^{\gamma_{n,m,s}-\frac{5}{2}n+\frac{3}{2}m+\frac{1}{2}s+1}\times
\\
&\times
\frac{(qz_1z_2)_{2n-1}(q^{3n+m+s}z_1z_2)_{s}(1-q^{6n+2m+2s-2}z_1z_2^2)}
{(q)_\infty (q)_{2n-2}(q)_{2n+m+2s}(q^{3n+m+s-1}z_2)_{n+s+1}
(q^{2n+m+2s+1}z_1z_2)_{n-1}(q^{4n+m+2s}z_1z_2^2)_\infty}\times
\\
&\times
[q^{\al_{n,m,s}+m}(z_1z_2^2)^{n+s}
(1:q^{4n+m+2s-1}z_2:q^{2n+m+2s}z_2)], \qquad n\ge 1,m\ge 1,s\ge 0,
\end{align*}

\begin{align*}
&D^{3n}A^{m}D(B+D+E)L^{s-1}v_\infty=(-1)^{n+m+s+1}(z_1z_2^3)^{n+s}z_2^{-1}
q^{\gamma_{n,m,s}-\frac{7}{2}n-\frac{3}{2}m-\frac{5}{2}s}\times
\\
&\times
\frac{(qz_1z_2)_{2n}(q^{3n+m+s+1}z_1z_2)_{s-1}(1-q^{6n+2m+2s}z_1z_2^2)}
{(q)_\infty(q)_{2n}(q)_{2n+m+2s-1}(q^{3n+m+s}z_2)_{n+s}(q^{2n+m+2s}z_1z_2)_{n}(q^{4n+m+2s+1}z_1z_2^2)_\infty}\times
\\
&\times
[q^{\al_{n,m,s}}(z_1z_2^2)^{n+s}(1:q^{-2n-m-2s}(z_1z_2)^{-1}:q^{2n}z_2)],
\qquad \quad n\ge 0,m\ge 1,s\ge 1,
\end{align*}

\begin{align*}
&D^{3n-2}CA^{m-1}D(B+D+E)L^{s}v_\infty=
(-1)^{n+m+s}(z_1z_2^3)^{n+s}
q^{\gamma_{n,m,s}-\frac{5}{2}n+\frac{1}{2}m-\frac{3}{2}s}(q)_\infty^{-1}\times
\\
&\times
\frac{(qz_1z_2)_{2n-1}(q^{3n+m+s}z_1z_2)_{s}(1-q^{6n+2m+2s-2}z_1z_2^2)}
{(q)_{2n-1}(q)_{2n+m+2s-1}(q^{3n+m+s-1}z_2)_{n+s+1}
(q^{2n+m+2s+1}z_1z_2)_{n-1}(q^{4n+m+2s}z_1z_2^2)_\infty}\times
\\
&\times
[q^{\al_{n,m,s}+m}(z_1z_2^2)^{n+s}
(1:q^{4n+m+2s-1}z_2:q^{2n-1}z_2)],\qquad n\ge 1,m\ge 1,s\ge 0,
\end{align*}

\begin{align*}
&D^{3n-1}CA^{m}D(B+D+E)L^{s-1}v_\infty=
(-1)^{n+m+s+1}(z_1z_2^3)^{n+s}z_2^{-1}
q^{\gamma_{n,m,s}-\frac{11}{2}n-\frac{3}{2}m-\frac{5}{2}s}\times
\\
&\times
\frac{(qz_1z_2)_{2n}(q^{3n+m+s+1}z_1z_2)_{s-1}(1-q^{6n+2m+2s}z_1z_2^2)}
{(q)_\infty(q)_{2n-1}(q)_{2n+m+2s-1}(q^{3n+m+s}z_2)_{n+s}
(q^{2n+m+2s}z_1z_2)_{n}(q^{4n+m+2s}z_1z_2^2)_\infty} \times
\\
&\times
[q^{\alpha_{n,m,s}}(z_1z_2^2)^{n+s}
(1:q^{-2n-m-2s}(z_1z_2)^{-1}:q^{-4n-m-2s}(z_1z_2)^{-1})], \;\;
n\ge 1,m\ge 0,s\ge 1,
\end{align*}

\begin{align*}
&D^{3n-3}CA^{m}D(B+D+E)L^{s}v_\infty=
(-1)^{n+m+s+1}(z_1z_2^3)^{n+s}z_2^{-1}
q^{\gamma_{n,m,s}-\frac{13}{2}n+\frac{1}{2}m-\frac{3}{2}s+2}\times
\\
&\times
\frac{(qz_1z_2)_{2n-2}(q^{3n+m+s}z_1z_2)_{s}(1-q^{6n+2m+2s-2}z_1z_2^2)}
{(q)_\infty(q)_{2n-2}(q)_{2n+m+2s}(q^{3n+m+s-1}z_2)_{n+s}
(q^{2n+m+2s+1}z_1z_2)_{n-1}(q^{4n+m+2s}z_1z_2^2)_\infty}\times
\\
&\times
[q^{\al+{n,m,s}+m}(z_1z_2^2)^{n+s}
(1:q^{-2n+1}(z_1z_2)^{-1}:q^{2n+m+2s}z_2)],\qquad n\ge 1,m\ge
0,s\ge 0,
\end{align*}

\begin{align*}
&D^{3n-2}E^{2m+1}L^{s}v_\infty=
(-1)^{n+m+1}(z_1z_2^2)^{n+s}z_2^m
q^{\delta_{n,m,s}-\frac{3}{2}n-\frac{1}{2}m}(q)_s^{-1}\times
\\
&\times
\frac{(q^{-2s}z_1)_s(qz_1z_2)_{n-1}}{(q)_{n-1}(q)_{2n+m+2s-1}
(q^{1-n-m-2s}z_1)_{n+m+2s-1}(q^{n+m+2s+1}z_2)_\infty
(q^{2n+m+2s}z_1z_2)_\infty}\times
\\
&\times
[q^{\beta_{n,m,s}}(z_1z_2^2)^{n+s}z_2^m
(1:q^{-2n-m-2s}(z_1z_2)^{-1}:q^{-n}z_1^{-1})],\qquad n\ge 1,m\ge 0,s\ge 0,
\end{align*}

\begin{align*}
&D^{3n-1}E^{2m+1}L^{s}v_\infty=
(-1)^{n+m}(z_1z_2^2)^{n+s}z_2^m
q^{\delta_{n,m,s}+\frac{1}{2}n+\frac{1}{2}m+\frac{1}{2}s}\times
\\
&\times
\frac{(q^{-2s}z_1)_s(qz_1z_2)_{n}}
{(q)_s(q)_{n-1}(q)_{2n+m+2s}(q^{-n-m-2s}z_1)_{n+m+2s}
(q^{n+m+2s+1}z_2)_\infty(q^{2n+m+2s+1}z_1z_2)_\infty}\times
\\
&\times
[q^{\beta_{n,m,s}}(z_1z_2^2)^{n+s}z_2^m
(1:q^{n+m+2s}z_1^{-1}:q^{2n+m+2s}z_2)],\qquad n\ge 1,m\ge 0,s\ge 0,
\end{align*}

\begin{align*}
&D^{3n}E^{2m-1}L^{s}v_\infty=
(-1)^{n+m}(z_1z_2^2)^{n+s}z_2^m
q^{\delta_{n,m,s}-\frac{1}{2}n-\frac{1}{2}m}\times
\\
&\times
\frac{(q^{-2s}z_1)_s(qz_1z_2)_{n}}
{(q)_s(q)_{n}(q)_{2n+m+2s-1}(q^{-n-m-2s+1}z_1)_{n+m+2s-1}
(q^{n+m+2s}z_2)_\infty(q^{2n+m+2s+1}z_1z_2)_\infty}
\\
&\times
[q^{\beta_{n,m,s}}(z_1z_2^2)^{n+s}z_2^m
(1:q^{n}:q^{-n-m-2s})],\qquad \qquad n\ge 0,m\ge 1,s\ge 0,
\end{align*}

\begin{align*}
&D^{3n-2}E^{2m+2}L^{s}v_\infty=
(-1)^{n+m}(z_1z_2^2)^{n+s}z_2^m
q^{\delta_{n,m,s}-\frac{3}{2}n-\frac{1}{2}m}\times
\\
&\times
\frac{(q^{-2s}z_1)_s(qz_1z_2)_{n-1}}
{(q)_s(q)_{n-1}(q)_{2n+m+2s-1}(q^{-n-m-2s}z_1)_{n+m+2s}
(q^{n+m+2s+1}z_2)_\infty(q^{2n+m+2s+1}z_1z_2)_\infty}\times
\\
&\times
[q^{\beta_{n,m,s}}(z_1z_2^2)^{n+s}z_2^m
(1:q^{n+m+2s}z_1^{-1}:q^{-n}z_1^{-1})],\qquad\qquad n\ge 1,m\ge 0,s\ge 0,
\end{align*}

\begin{align*}
&D^{3n-1}E^{2m}L^{s}v_\infty=
(-1)^{n+m}(z_1z_2^2)^{n+s}z_2^m
q^{\delta_{n,m,s}-\frac{3}{2}n-\frac{1}{2}m}\times
\\
&\times
\frac{(q^{-2s}z_1)_s(qz_1z_2)_{n}}
{(q)_s(q)_{n-1}(q)_{2n+m+2s-1}(q^{-n-m-2s+1}z_1)_{n+m+2s-1}
(q^{2n+m+2s}z_1z_2)_\infty(q^{n+m+2s}z_2)_\infty}\times
\\
&\times
[q^{\beta_{n,m,s}}(z_1z_2^2)^{n+s}z_2^m
(1:q^{-2n-m-2s}(z_1z_2)^{-1}:q^{-n-m-2s})], \qquad n\ge 1,m\ge 1,s\ge 0,
\end{align*}

\begin{align*}
&D^{3n}E^{2m}L^{s}v_\infty=
(-1)^{n+m}(z_1z_2^2)^{n+s}z_2^m
q^{\delta_{n,m,s}+\frac{3}{2}n+\frac{1}{2}m+2s}\times
\\
&\times
\frac{(q^{-2s}z_1)_s(qz_1z_2)_{n}}
{(q)_s(q)_{n}(q)_{2n+m+2s}(q^{-n-m-2s+1}z_1)_{n+m+2s-1}
(q^{n+m+2s+1}z_2)_\infty(q^{2n+m+2s+1}z_1z_2)_\infty}\times
\\
&\times
[q^{\beta_{n,m,s}}(z_1z_2^2)^{n+s}z_2^m
(1:q^{n}:q^{2n+m+2s}z_2)],\qquad\qquad n\ge 0,m\ge 1,s\ge 0.
\end{align*}
Our convention 
for expansions of RHS
of these formulas into power series can be now written in the form:
\be
|q|\ll 1,\qquad |z_1^{-1}|\ll 1,\qquad |z_2|\ll 1, \qquad |z_1z_2|\sim 1.
\ee


\end{document}